
\mathsurround=2pt
\def\fn{$f:({\bf C}^n , 0) \rightarrow ({\bf C}, 0)$}
\def\gn{$g:({\bf C}^n , 0) \rightarrow ({\bf C}, 0)$}
\def\Cn{{\bf C}^n}
\def\[f]{\lbrack f\rbrack}
\def\ADE{$A_k , D_k, E_6 ,  E_7 , E_8$}
\magnification=\magstep1
\baselineskip=0.7cm
\hsize=12cm

\centerline {\bf ON THE A-D-E CLASSIFICATION OF}
\centerline {\bf THE SIMPLE SINGULARITIES OF FUNCTIONS}

\bigskip
\centerline {MIKHAIL ENTOV}

\bigskip
\bigskip
\centerline {\bf Introduction.}

\bigskip

     One of the most marvellous results of singularity theory is
the classification  of  the  simple  isolated  singularities  of
analytic  functions by the Weyl groups (and the irreducible root
systems) of the types $A_k, D_k, E_6, E_7, E_8$.

     There are various interesting relations between the  simple
singularities and the algebraic objects (simple Lie groups, root
systems,  Coxeter groups) of the types $A_k, D_k, E_6, E_7, E_8$
(see, e.g., [Br1], [Du], [Gr], [Sl]). All these relations deeply
involve the normal forms of the simple singularities obtained by
V.I.Arnold in his paper [Ar1].

     In this paper we consider one of such relations  not  using
the  normal  forms. We dwell on the construction that associates
to a singularity its monodromy group (for a suspension  with  an
odd number of variables).

     It  is well known that this con\-struc\-tion as\-so\-ciates
to each of the simple sin\-gu\-la\-ri\-ties of the types $A_k  ,
D_k,  E_6,  E_7,  E_8$  the  Weyl  group  of  the same type, and
moreover, the monodromy group of the  singularity  is a finite
group  generated  by  reflections  only  if the singularity is
simple. The normal forms are essentially used in the  proofs  of
these facts (see [Ar1], [Ar2],  [Ar4],  [AGV1],  [AGV2],[Tju1]),
and  the  coincidence  of  the  classifications  of  the  simple
singularities and of the (irreducible) Weyl groups  looks  just
like  the  mysterious  coincidence  of the lists obtained by the
independently proved classification theorems.

     That gave rise  to  the  natural  problem:  to  reduce  the
classification  of  the simple singularities directly, not using
the normal forms, to the classification of the irreducible  Weyl
groups. This problem was repeatedly mentioned by V.I.Arnold (see
e.g. [Ar1], [Ar2], [Ar6]).

     In  this paper we show how this problem can be solved. More
exactly, we prove (see theorems I and II), not using the  normal
forms,  that

\medskip

{\sl
\noindent
a) the monodromy group of a
simple singularity is a Weyl group;

\noindent
b) if the monodromy group of a singularity is finite, then it is
an irreducible Weyl group of one of the types  $A_k,  D_k,  E_6,
E_7, E_8$, and the singularity is simple;

\noindent
c) if two simple singularities have isomorphic monodromy groups,
then they are equivalent.}

\medskip

\noindent
     We reduce the proofs of a) and  b)  to  the  proof  of  the
assertion  that  the  simple sin\-gu\-la\-ri\-ties coincide with
the elliptic ones, i.e. with  ones  with  definite  intersection
form  on  the  homologies  of the Milnor fiber (for a suspension
with an odd number of variables), and we obtain this  assertion
proving  that  a singularity is simple (elliptic) if and only if
the mixed Hodge  structure  in  its  vanishing  cohomologies  is
trivial  or,  in  other words, the length of the spectrum of the
singularity is less than one.

\smallskip

     In the course of proving of a) and b) we also find that

\medskip

\noindent
{\sl a simple singularity is stably  equivalent  to
the  singularity  of  a   quasihomogeneous   function   of   two
variables.}

\medskip

     We  show (not using the normal forms) that

\medskip

\noindent
{\sl the monodromy operator of a simple singularity is a Coxeter
element of the corresponding Weyl group.}

\medskip

     From  this,  by  virtue  of  the purely algebraic result of
P.Deligne obtained in [De], we deduce  that

\medskip

\noindent
{\sl the Dynkin diagram of a simple singularity is the canonical
Dynkin   diagram  of  the  corresponding  Weyl  group  for  some
distinguished basis.}

\medskip

     Then one can easily find the normal  forms  of  the  simple
singularities  and determine the types of the corresponding Weyl
groups as it is described in Appendix 1.

\smallskip

Our approach is that instead of the use of the normal  forms  we
apply   many   other   general  and  very  powerful  results  of
singularity theory and of the  theory  of  Weyl  groups  to  the
particular  case  of a simple (an elliptic) singularity. Most of
these results were obtained  after  V.I.Arnold  had  found  the
normal forms of the simple singularities.

     We also  use  the  new  result  (see  theorem  1)  that

\medskip

\noindent
{\sl  under the canonical identification of the local algebra of
a singularity with the tangent space to the base of a miniversal
deformation of  the  singularity  at  zero,  the  class  of  the
singularity  in  its local algebra always belongs to the tangent
cone to the stratum \hbox{$\mu {\hbox{\rm =const}}$}.}

\medskip

     This theorem has been proved  jointly  by  J.H.M.Steenbrink
and  the  author.  For the case of the singularity of a function
which is nondegenerate with respect to its  Newton  diagram  the
result  can  be  sharpened  (see  Appendix  2)  in  the way that

\medskip

\noindent
{\sl   under    the    mentioned    above
identification,  the  class  of  the  singularity  in  its local
algebra is the tangent vector  to  a  {\it  linear}  1-parameter
\hbox{$\mu    {\hbox{\rm    =const}}$}    deformation   of   the
singularity.}

\medskip

     The proof of the assertion that the monodromy operator of a
simple singularity is a Coxeter  element  of  the  corresponding
Weyl  group  is  based  in  particular on the interesting purely
algebraic remark (see proposition 2.1) that

\medskip

\noindent
{\sl  the  Coxeter elements of a Weyl group (of one of the types
\ADE ) are the only elements of maximal length (i.e. those  that
can  be  written as an irreducible product of the maximal number
of reflections) with trace equal to $-1$.}

\medskip

     I am deeply grateful to J.H.M.Steenbrink for the help,  the
numerous  fruitful  discussions and many useful comments on this
text. The proof of theorem 1 given  in  this  paper  belongs  to
J.H.M.Steenbrink.  Also  the  proof of the implication ${\rm 2)}
\Rightarrow {\rm 6)}$ (in theorem I) has  been  inspired  by  my
conversation  with  J.H.M.Steen\-brink.  I  am  also grateful to
V.I.Ar\-nold,   A.A.Bei\-lin\-son,    R.V.Bez\-ru\-kav\-ni\-kov,
Yu.G.Makh\-lin,       A.G.Ho\-van\-sky,       V.A.Vas\-sil\-yev,
M.S.Ver\-bit\-sky and especially to E.J.N.Lo\-oijen\-ga for  the
useful  discussions  and for helping me in various ways. I thank
E.Brieskorn for many valuable remarks  he  made  reviewing  this
paper.  The  first proof of theorem II
in \S 7  was communicated to me by
E.Brieskorn.

     This paper has been completed during  my  participation  in
the 1992--93 Masterclass  on Geo\-met\-ry organized by the University of
Utrecht and the Ma\-the\-ma\-ti\-cal Research Institute  in  The
Netherlands. I thank these institutions for the hospitality.

\bigskip
\bigskip

     \centerline {\bf \S 1. Terminology.}

\bigskip

In  this  paragraph  we  introduce the notation and recapitulate
some constructions of singularity theory and the theory of  Weyl
groups.

\smallskip

\noindent
{\bf 1.1. Some notions and constructions of singularity theory.}
Basic  notions  of  singularity  theory  can be found in [AGV1],
[AGV2], [AGLV].

     In this paper we deal only with isolated  singularities  of
analytic functions of several complex variables.

     By $O_n$ we denote the space of germs of analytic functions
of $n$ comp\-lex variables at the point ${0\in \Cn}$. Let \fn be
an  isolated  singularity  of  an analytic function at the point
${0\in \Cn}$. By $I_f$ we denote the gradient ideal of $f$, i.e.
the ideal in $O_n$ generated by all partial derivatives of  $f$
at  zero.  By  $Q_f$  we  denote the local algebra of $f$ , i.e.
$Q_f=O_n/I_f$. By $\[f]$ we denote the class of a  germ  $f$  in
its local algebra $Q_f$.

     Now we fix such sufficiently small neighbourhoods  of  zero
\break  $U={\lbrace  z\mid  \  {\parallel  z  \parallel } < \rho
\rbrace} \subset \Cn$ and $T= {{\lbrace t \mid \ }{{\parallel  t
\parallel} < \delta \rbrace }\subset {\bf C}}$ (see [AGV2]) that

\medskip

\noindent
a) the point $0\in \Cn$ is the  only  critical
point of $f$ in the ball $U$;

\noindent
b)  the  hypersurface $f^{-1} (t)$ is nonsingular inside $U$ and
intersects transversally the  boundary  of  $U$  for  all  $t\in
T\setminus 0$.

\medskip

\noindent
Let  $X_t=f^{-1}  (t)  \cap  U$,  $t\in  T\setminus 0$, denote a
nonsingular level of $f$ near the critical point  $0$.  We  also
fix  a  point $t\in T\setminus 0$. The manifold $X_\ast=X_t$ has
the  homotopy  type  of  the   bouquet   of   $\mu$   $(n-1){\rm
-dimensional}$  spheres,  where $\mu$ is the Milnor number (or
the multiplicity) of the singularity $f$  (see  [Mi]).  So  $H_k
(X_\ast)=0$   if   $k\not=  0,  n-1$  and  $H_{n-1}(X_\ast) \cong {\bf
Z}^\mu$ (for $n=1$  here and further
one should  consider the reduced homology
group ${\widetilde H}_0 (X_\ast)$).

     We  take  such  a small perturbation $f_\varepsilon$ of the
singularity $f$ that for all  sufficiently  small  $\varepsilon$
the  function  $f_\varepsilon$  has exactly $\mu$ Morse critical
points inside $U$ with different  critical  values  inside  $T$.
Let's   fix   any   such   $\varepsilon$.  A  nonsingular  level
${f_\varepsilon}^{-1}(t)\cap U$ is  homeomorphic  to  the  level
$X_\ast$  of the function $f$. We take a noncritical value $t_0$
of $f_\varepsilon$ on the boundary of $T$ and then choose  $\mu$
non-self-intersecting  and  mutually  nonintersecting (except at
the point $t_0$) paths ${\gamma}_1,\ldots  ,{\gamma}_\mu$  going
inside  $T$  from  the  point  $t_0$  to  the critical values of
$f_\varepsilon$. The paths are numbered in  the  same  order  in
which they emanate from the point $t_0$ counting clockwise. Such
a  system  of  paths  is  called  a {\it distinguished system of
paths} (see [Ga1]  or[AGV2]).  The  cycles  ${\Delta}_1,  \ldots
,{\Delta}_\mu$  vanishing  along  the  paths ${\gamma}_1,\ldots
,{\gamma}_\mu$ form a  basis  (a  so-called  {\it  distinguished
basis})  of  $H_{n-1}  (f_{\varepsilon}^{-1} (t_0) \cap U) \cong
H_{n-1}(X_\ast)$ (see [Br2], [Ga1], [Lam] or [AGV2]).

     We assume that the number of variables  $n\equiv  3\  ({\rm
mod}\  4)$.  Otherwise one takes a suitable suspension of $f$ --
for all such suspensions the monodromy groups are isomorphic and
the quadratic forms are the same (see [AGV2]), so

\medskip

\noindent
{\it further, whenever we speak about the monodromy group or the
quadratic form of a singularity, we mean  that  the  appropriate
(as it is mentioned above) su\-spen\-sion is considered.}

\medskip

\noindent
Then (see [AGV2]) the intersection form in  the  homology  group
\break
$H_{n-1}(X_\ast, {\bf  R}) \cong {\bf  R}^\mu$  is  symmetric  and  it
defines  a  quadratic form called {\it the quadratic form of the
singularity}  $f$.  The   index   of   self-intersection   of   a
va\-ni\-shing  cycle  is  equal to $-2$, and the monodromy group
$\Gamma$  of  the  singularity   $f$   is   generated   by   the
Picard-Lefschetz  transformations $h_i$ related to the vanishing
cycles ${\Delta}_i$ (see  [Gu]  or  [AGV2]):  $$  h_i  :  \sigma
\rightarrow \sigma + (\sigma \circ {\Delta}_i) {\Delta}_i,\leqno
(1)$$ where $\sigma \in H_{n-1}(X_\ast)$, $i=1,\ldots ,\mu$, and
$(\cdot   \circ  \cdot)$  denotes  the  intersection  form.  The
monodromy operator of the singularity $f$ is  the  product  $h_1
\cdot \ldots \cdot h_\mu$.

\smallskip

\noindent
{\bf Definition} (see [Ar2]). A singularity is said to  be  {\it
elliptic} if its quadratic form is negative definite.

\smallskip

\noindent
{\bf  Remark.} This notion of an elliptic singularity should not
be confused with one introduced in  the  works  of  K.Saito  and
E.J.N.Looijenga.

\smallskip

     The  braid group $Br(\mu )$ acts transitively on the set of
the  systems  of  distinguished  paths  (considered  up   to   a
homotopy).  Namely, if $b_1,\ldots ,b_{\mu -1}$ are the standard
generators  of  $Br(\mu)$  then  $b_i$  transfers   the   system
${\gamma}_1,\ldots      ,{\gamma}_\mu$     to     the     system
${\gamma}_1^{\prime},\ldots ,{\gamma}_\mu^{\prime}$: \noindent $
{\gamma}_j^{\prime}   =   {\gamma}_j   $,   $j\ne   i,   {i+1}$,
${\gamma}_i^{\prime}         =         {\gamma}_{i+1}$,        $
{\gamma}_{i+1}^{\prime}$  is  homotopic  to   ${\gamma}_i   \cup
{\varrho}_{i+1}^{-1}$,  \noindent  where  ${\varrho}_i$  is  the
simple loop going from the fixed  point  $t_0$  along  the  path
${\gamma}_i$  to a point nearby the end of
${\gamma}_i$,  then  once  counterclockwise   around   the   end   of
${\gamma}_i$  and  then  back  along  ${\gamma}_i$  to $t_0$. By
formula (1) this action of $Br(\mu )$  provides  the  action  of
$Br(\mu)$ on the set of the distinguished bases (see [Lo1], [Gu]
or [AGV2]):


$$\leqalignno{    \hbox{$b_i    :$}&\    ({\Delta}_1,     \ldots
,{\Delta}_\mu  )  \mapsto  &(2)\cr \hfil&\ ({\Delta}_1, \ldots ,
{\Delta}_{i-1}, {\Delta}_{i+1}, {\Delta}_i -  ({\Delta}_i  \circ
{\Delta}_{i+1})     {\Delta}_{i+1},    {\Delta}_{i+2},    \ldots
,{\Delta}_\mu).\cr} $$

\medskip

\noindent
{\bf  1.2.  Some  notions  on  Weyl  groups  and  finite  groups
generated  by reflections. } By a {\it finite group generated by
reflections} we mean a finite group generated by reflections  in
a  Euclidean  vector  space. If $\Re$ is a (reduced) root system
then by $W(\Re )$ we denote the corresponding Weyl group.  Basic
notions  on  groups generated by reflections and Weyl groups can
be found in [Bou].

     Now let $W(\Re )$ be a Weyl group (of one of the types \ADE
) of rank $\mu$ and let $S$ denote the  set  of  all  \hbox{$\mu
{\hbox{\rm -tuples}}$} $(s_1,\ldots , s_\mu )$ of reflections in
$W(\Re )$ such that

\smallskip

\noindent
i) $s_1,\ldots , s_\mu $ generate $W(\Re )$;

\noindent
ii)  the  roots  corresponding  to  $s_1,\ldots  ,  s_\mu  $ are
linearly independent and span $\Re $ over ${\bf Z}$.

\smallskip

\noindent
Then the braid group $Br(\mu )$ acts on $S$ in a way similar  to
(2)  (see [Lo1]) : $$\leqalignno{ \hbox{$b_i :$} &\ (s_1, \ldots
,s_{i-1} ,s_i ,s_{i+1} ,s_{i+2} ,\ldots ,s_\mu ) \mapsto &(3)\cr
\hfil &\ (s_1, \ldots ,s_{i-1}  ,s_{i+1}  ,s_{i+1}  s_i  s_{i+1}
,s_{i+2} ,\ldots ,s_\mu ). \cr }$$

\noindent
One  can  easily notice that this action leaves the product $s_1
\cdot \ldots \cdot s_\mu $ invariant.

\vfil
\eject

\bigskip
\bigskip

\centerline {\bf \S 2. Statement of Results.}

\bigskip

\noindent
We use the notation and the agreements introduced in \S 1.

\noindent
The main results of this paper are the following two theorems.

\smallskip

\noindent
{\bf  Theorem  I.}  {\it  For  any  singularity  the   following
conditions are equivalent:

\smallskip

\noindent
1) The singularity is simple;

\noindent
2) The singularity is elliptic;

\noindent
3) The monodromy group of the singularity is finite;

\noindent
4) The monodromy group of the singularity is isomorphic
to a Weyl group of one of the types \ADE;

\noindent
5)  The  mixed  Hodge structure in the vanishing cohomologies of
the singularity is trivial (i.e. both the Hodge and  the  weight
filtrations do not contain any nontrivial subspaces);

\noindent
6) The length of the spectrum  of  the  singularity  is
less than one.}

\medskip

\noindent
{\bf   Theorem  II.}  {\it  If  two  simple  singularities  have
isomorphic monodromy groups then they are (stably) equivalent.}

\smallskip

     The normal forms of the simple singularities are  not  used
in  the  proofs  of theorems I,II. In the course of the proof of
theorem I we also obtain

\smallskip

\noindent
{\bf Corollary 1.}

\noindent
{\it  A simple (an elliptic) singularity is stably equivalent to
the  singularity  of  a   quasihomogeneous   function   of   two
variables.}

\smallskip

\noindent
{\bf Remark.} The equivalence  ${\rm  1)}  \Leftrightarrow  {\rm
6)}$ was conjectured by K.Saito in [Sa3].

\medskip

\noindent
To  prove theorem I we use the following result obtained in this
paper (see \S 4).

\smallskip

\noindent
{\bf Theorem 1.} {\it Let $F:  \Cn  \times  \Lambda  \rightarrow
{\bf  C}  $  be  a  miniversal  deformation of $f$, and let $T_0
(\Lambda)$ denote the tangent space to $\Lambda$ at  zero.  Then
if  $v  \in  T_0 (\Lambda )$ is mapped onto $\[f] \in Q_f$ under
the canonical identification of $ T_0 (\Lambda  )$  with  $Q_f$
then  $v$  lies  in  the  tangent cone to the stratum \hbox{$\mu
{\hbox{\rm =const}}$} in the space $\Lambda$. }

\smallskip

This  result  fits  with  the  theorem  of   A.N.Varchenko   and
S.V.Chmutov  on  the  tangent  cone  to  the  stratum \hbox{$\mu
{\hbox{\rm =const}}$} of a singularity (see [VC]). For the proof
of theorem 1 see \S 4. Moreover, the following fact is also true
(see Appendix 2).

\smallskip

\noindent
{\bf Theorem 1'.} {\it Let's fix a monomial basis of $Q_f$  over
{\bf C} and consider $\[f]$ as a linear combination of the basic
monomials. If $f$ is nondegenerate with respect  to  its  Newton
diagram (see {\rm [Kou]} or {\rm [AGV2]}), then for small $t$ the
linear  deformation  $f_t  =  f + t{\[f]}$, $t\in {\bf C}$, is a
deformation with constant multiplicity.}

\smallskip

\noindent
We also prove (not using the normal forms) that

\noindent
{\bf Theorem 2.} {\it The monodromy operator  of  a  simple  (an
elliptic)  singularity with the monodromy group $W$ is a Coxeter
element of the Weyl group $W$.}

\smallskip

\noindent
{\bf Corollary 2.} {\it  Let  $f$  be  a  simple  (an  elliptic)
singularity. Let a Weyl group $W$ be the monodromy group of $f$.
Then  there  exists  such  a distinguished basis that the Dynkin
diagram of $f$ with respect to
this basis is the canonical Dynkin diagram  of
the Weyl group $W$.}

\smallskip

\noindent
To prove theorem 2 we use the following algebraic result.


\smallskip

\noindent
{\bf Proposition 2.1.} {\it Let
$W(\Re)$ be an irreducible Weyl group of the rank $\mu$ where $\Re$ is
a (reduced) root system of one of the types \ADE .  Let  $s=s_1  \cdot
\ldots  s_\mu$  be  a  product  of  $\mu$  reflections   in   $W(\Re)$
corresponding to some linearly independent roots that span  the  root
system $\Re$ over ${\bf Z}$ (i.e. $(s_1, \ldots ,s_\mu) \in S$ following
the notation from \S 1).
Then $s$ is a Coxeter element  of  $W(\Re
)$ if and only if its trace is equal to $-1$.}

\smallskip

\noindent
 Theorem 2 and corollary 1 also are necessary to  check  that  the
results we refer to while proving theorem II can  be  in  fact  proved
without the normal forms (see the remark in \S 7).

\smallskip

     Finally we show how using corollary 1 one  can  find  the  normal
forms of the simple singularities and then, by virtue  of  theorem  2,
determine the types of the corresponding Weyl groups (see Appendix 1).

\smallskip

     The scheme of the paper is as follows.  In  \S  3  we  prove  the
equivalence of conditions 2), 3) and 4) of theorem I. In \S 4 we prove
theorem 1. In \S 5 we prove the equivalence of conditions 1),  2),  5)
and 6) of theorem I and corollary 1. In \S 6  we  prove  proposition
2.1, theorem 2 and corollary 2. In \S 7 we prove \allowbreak theo\-rem
II.

\bigskip
\bigskip

\centerline  {\bf  \S  3.  The  Proof  of  the  Equivalence   of}
\centerline {\bf Conditions 2), 3) and 4) in Theorem I.}

\bigskip

We   shall   prove   the   implications   ${\rm    2)}
\Leftrightarrow {\rm 3)}$ and ${\rm 3)}  \Rightarrow  {\rm  4)}$.  The
implication ${\rm 4)} \Rightarrow {\rm 3)}$ is obvious. In  the  proof
we use the following well-known facts.

\medskip

\noindent
{\bf Proposition 3.1} (see  [Sa2]).  {\it  The
quadratic form of a singularity $f$ is the unique  integral  symmetric
even (i.e. with values in ${2}${\bf Z}) quadratic  form  on  ${\bf
Z}^\mu \cong H_{n-1} (X_\ast) $ invariant under the action of the monodromy
group $\Gamma$ and such that $-2$ is a value of it.}

\smallskip

\noindent
{\bf Proposition 3.2}  (see  [Gu],  [AGV2]).
{\it The monodromy group of a singularity $f$ acts transitively on the
set of the vanishing cycles of $f$,  i.e.  for  any  vanishing  cycles
$\Delta$  and  ${\Delta}^\prime$  there  exists  an  element  of   the
monodromy group taking $\Delta$ into $\pm {\Delta}^\prime$. }

\medskip

\noindent
{\bf The proof of the implications.}

\smallskip

\noindent
${\rm 2)}\Rightarrow {\rm 3)}$

\noindent
The monodromy group $\Gamma$
is a subgroup of the automorphism group of the
integral lattice
$H_{n-1}(X_\ast, {\bf Z}) \cong {\bf Z}^\mu \subset
H_{n-1}(X_\ast, {\bf R}) \cong {\bf R}^\mu$. The elements of $\Gamma$
also preseve the intersection form. Therefore if the intersection
form is definite then $\Gamma$ is a discrete subgroup of a compact
group and hence finite.

\smallskip

\noindent
${\rm 3)}\Rightarrow {\rm 2)}$

\noindent
Let  $j$  be  a  positive  definite  inner  product  on
$H_{n-1}(X_\ast, {\bf R}) \cong {\bf R}^\mu$ invariant under the  action  of
$\Gamma$ (such a  product  exists  because  $\Gamma$  is  finite).  By
proposition 3.2, all  vanishing  cycles  have  the  same  length  with
respect to $j$. Therefore normalizing $j$, if necessary, we can assume
that all vanishing cycles have the length $2$  with  respect  to  $j$.
Since one can choose a basis of  $H_{n-1}(X_\ast)$  from  the  set  of
vanishing cycles (see \S 1), the form $j$ restricted  to  the  lattice
$H_{n-1}(X_\ast, {\bf Z})\cong {\bf  Z}^\mu  \subset  H_{n-1}(X_\ast,  {\bf
R})\cong {\bf R}^\mu$ is integral, even and taking a value $2$.  Hence,  by
proposition 3.1, the quadratic  form  on  $H_{n-1}(X_\ast,  {\bf  R})$
defined by the bilinear form $-j  $  is  the  quadratic  form  of  the
singularity. So the singularity is elliptic.

\smallskip

\noindent
${\rm 3)}\Rightarrow {\rm 4)}$

\noindent
As it has been proved above, if the monodromy group  is
finite,  then  the  singularity  is   elliptic.   Hence
the Picard-Lefschetz transformations are the reflections in the Euclidean
vector space $H_{n-1}(X_\ast, {\bf R})$ and
the vanishing cycles form a root system,  so
the monodromy group is {\it a Weyl group}. It is irreducible because
of proposition 3.2. Also by proposition 3.2, all vanishing cycles  have
the same length, so the irreducible (reduced) root  system  they  form
can be only of one of the types \ADE. The implication is proved.

\bigskip
\bigskip

\centerline {\bf \S 4. The Proof of Theorem 1.}

\bigskip

     If $\[f] = 0$ then the statement is trivial. Let $\[f] \not=  0$.
Then one may choose $e_1 = 1, e_2,\ldots ,e_\mu =f$,
$({e_i
\in O_n}, i=1,\ldots ,\mu)$, mapping to a basis of $Q_f$ over {\bf C}.
Then $H: \Cn \times {\bf C}^\mu \rightarrow {\bf C}$, $H(z, \alpha)  =
f(z) + {\alpha}_1 e_1 (z) +\ldots + {\alpha}_\mu e_\mu  (z)  $,  $z\in
{\bf C}^n$, $\alpha =  ({\alpha}_1,\ldots  ,{\alpha}_\mu  )\in  {\bf
C}^\mu$, is a  miniversal  deformation  of  $f$.  So  there  exists  a
biholomorphism $g: \Lambda \rightarrow  {\bf  C}^\mu$  such  that  the
deformation $F$ is equivalent to the
one induced from $H$ by $g$. One  can
easily check that $dg_0 (v)
= (0,\ldots ,0,1)$, where $dg_0$ is the differential  of  $g$  at
zero. It is also clear that $g$ maps the stratum \hbox{$\mu {\hbox{\rm
=const}}$} in $\Lambda$ to the stratum\hbox{$\mu {\hbox{\rm =const}}$}
in ${\bf C}^\mu$. So it suffices to show that the vector  $  (0,\ldots
,0,1)\in  T_0  {\bf  C}^\mu$  lies  in  the  tangent   cone   to   the
stratum\hbox{$\mu {\hbox{\rm =const}}$}  in  ${\bf  C}^\mu$.  This  is
trivial because  the  family  $(1+t)f$,  $t\not=  -1$,  is  \hbox{$\mu
{\hbox{\rm -constant}}$}. The theorem is proved.

\bigskip
\bigskip

\centerline  {\bf  \S  5.  The  Proofs  of  the  Equivalence   of
                           Conditions }
\centerline {\bf 1), 2), 5)  and  6)  in  Theorem  I  and
Corollary 1.}

\bigskip

     Firstly we recall some known facts  of  singularity  theory  that
will be used in the proof.

\smallskip

\noindent {\bf Proposition 5.1} ("The Morse lemma with
parameters" -- see [Ar1] or [AGV1]). {\it  In  a  neighborhood  of  a
critical point of corank $k$ a holomorphic function
\break
\fn is  equivalent
to  a  function  ${\tilde  f}(z_1,\ldots  ,z_k)  +  z_{k+1}^2+\ldots
+z_n^2$, where the second differential of  ${\tilde  f}$  at  zero  is
equal to zero.}

\smallskip

\noindent
{\bf Proposition 5.2.} {\it  If  the  second
differential of the singularity \break \fn at zero is equal  to  zero,
then the positive (the negative) index of  inertia  of  the  quadratic
form of the singularity  $f$  is  not  less  than  the  positive  (the
negative) one of the quadratic form of the  singularity  ${\Theta  }_n
(z_1,\ldots ,z_n) = z_1^3 + \ldots +z_n^3$.}

\smallskip

\noindent
{\bf The  proof  of  proposition  5.2.}  The
proposition follows, for instance, from the results of G.N.Tjurina  --
see [Tju1], \S 1, proposition 2 and theorem 1.

\smallskip

\smallskip

\noindent
{\bf Proposition 5.3} (see [Ga2]). {\it  The
modality of the singularity is always one less than the  dimension  of
the  stratum \hbox{$\mu  {\hbox{\it  =const}}$}  in  the  base   of   a
miniversal deformation of the sin\-gu\-la\-ri\-ty.}

\smallskip

\noindent
{\bf Proposition 5.4} (see [V5]).  {\it  The
codimension of the stratum \hbox{$\mu  {\hbox{\it  =const}}$}  in  the
base of a miniversal deformation of a singularity $f$ is not less than
the number of spectral numbers of $f$ that are  less  than  $l_1  +1$,
where $l_1$ is the minimal spectral number of $f$.}


\smallskip

\noindent
{\bf Proposition  5.5}  (see  [St1],  [V4]).
{\it Let's assume that the intersection form of the singularity $f$ is
nondegenerate. Let $({\mu}_{+}, {\mu}_{-})$ denote its signature. Then

\noindent
a) All spectral numbers of $f$ are not integer;

\noindent
b) The index ${\mu}_{+}$ (${\mu}_{-}$) is equal to  the
number of the spectral numbers of $f$ with odd (even) integral part.}

\smallskip

\noindent
{\bf Proposition 5.6} (see [V2]).  {\it  If
$\lbrace l_i \rbrace$, $ i=1,\ldots ,\mu $, is  the  spectrum  of  the
singularity \fn , then $\lbrace l_i + 1/2 \rbrace$,  $i=1,\ldots  ,\mu
$, is the spectrum of the singularity $f + z_{n+1}^2 :({\bf C}^{n + 1}
, 0) \rightarrow ({\bf C}, 0)$.}

\smallskip

\noindent
{\bf Proposition 5.7} (see [V3]).  {\it  Let
$\lbrace f \rbrace $ denote the operator of multiplication by  $f$  in
the local algebra $Q_f$. Then if a number  $j$  is  greater  than  the
length of the spectrum of $f$, then ${\lbrace f\rbrace }^j =0$.}


\smallskip

\noindent
{\bf  Proposition  5.8}  (see
[Sa1]). {\it A singularity $f$ is equivalent to the singularity  of  a
quasihomogeneous function if and only if the class $\[f]  =0$  in  the
local algebra $Q_f$.}

\smallskip

We   recall   that   for   the   singularity   of   a
quasihomogeneous function \break {\fn} of degree one a monomial from a
monomial  basis  of  $Q_f$  over  {\bf  C}  is  called   {\it   upper}
(respectively, {\it diagonal} or {\it lower}) if its  quasihomogeneous
degree is greater than one (respectively, equal to one  or  less  than
one). The total number of upper (diagonal, lower) basic monomials  in
a monomial basis does not depend on the choice  of  such  a  basis  of
$Q_f$ (see [Ar3], [AGV1]).

\smallskip

\vfil
\eject

\noindent
{\bf Proposition 5.9} (see [Ar3],[V5]). {\it
The modality of the singularity of a quasihomogeneous function $f$  is
equal to the total number of upper and diagonal basic monomials of the
local algebra $Q_f$.}

\smallskip

\noindent
{\bf Proposition 5.10} (see [St2]). {\it Let
\fn be the singularity of a quasihomogeneous function  of  degree  $1$
with the  weights  $\nu  =  ({\nu}_1,\ldots  ,{\nu}_n  )$.  Let  ${\bf
z}^{{\bf k}_i}$, $i=1,\ldots , \mu$, be a monomial basis of the  local
algebra $Q_f$ over {\bf C}. Then the  spectrum  of  $f$  is  the  set
$\lbrace \langle {\bf k}_i +{\bf 1}, \nu \rangle - 1\rbrace  $,  $i=1,\ldots
,\mu $,where ${\bf 1}=(1,\ldots ,1)$.}

\medskip

\noindent
{\bf The proof of the equivalences
${\bf  {\rm 1)} \Leftrightarrow {\rm 2)} \Leftrightarrow {\rm 5)}
\Leftrightarrow
{\rm 6)}}$ in theorem I.}

\smallskip

\noindent
${\rm  5)} \Leftrightarrow {\rm 6)}$

\smallskip

\noindent
Follows  from
the definitions and the symmetries of the spectrum (see [St1], [V2]).

\medskip

\noindent
${\rm  1)}  \Rightarrow  {\rm  6)}$

\smallskip

\noindent
If  the  singularity  $f$  is  simple,  then,  by
proposition 5.3, the stratum \hbox{$\mu {\hbox{\rm  =const}}$}  in  the
base of a miniversal deformation of $f$ is 1-dimensional and  consists
only of the coordinate axis ${\lambda}_1$. Hence, by  theorem  1,  the
class $\[f] =0$ and,  by  proposition  5.8,  the  singularity  $f$  is
equivalent to the singularity of a quasihomogeneous function.  So  let
\fn be a simple quasihomogeneous singularity of degree  $1$  with  the
weights $\nu  =  ({\nu}_1,\ldots  ,{\nu}_n  )$.  Then,  by  virtue  of
proposition 5.10, it can  be  easily  seen  that  the  length  of  the
spectrum of $f$ is equal to the  maximal  quasihomogeneous  degree  of
basic monomials of $Q_f$. Therefore, by  proposition  5.9,  one  finds
that the length of the spectrum of $f$ is less
than $1$, i.e. all spectral numbers lie  in  the  interval  $(n/2
-3/2 , n/2 - 1/2)$.

\medskip

\noindent
${\rm  6)}  \Rightarrow  {\rm  1)}$

\smallskip

\noindent
Follows from propositions 5.4 and 5.3.

\medskip

\noindent
${\rm 2)} \Rightarrow {\rm  6)}$  (cf.  [Tju1])

\smallskip

\noindent
Let \fn be an elliptic  singularity.  We
claim that $f$ is stably equivalent to the singularity of  a  function
of two variables.

     Indeed, let $k$ be the corank of $f$. Then, by  proposition  5.1,
the singularity $f$ is stably  equivalent  to  the  singularity  of  a
function ${\tilde f}:({\bf C}^k , 0) \rightarrow ({\bf C},  0)$,  such
that the second differential of ${\tilde f}$ at zero is equal to zero.
Using propositions 5.10 and 5.5  one  computes  the  spectrum  of  the
singularity ${\Theta}_k $ and sees that the singularity $ {\Theta}_k $
is elliptic if and only if $k \leq 2$. Therefore, by proposition  5.2,
$f$ is equivalent to the singularity of a function ${\tilde f}$ of two
variables.

     By the symmetry  and by proposition 5.5, we find  that
the spectrum of the singularity ${\tilde f} +  z_3  ^2$  lies  in  the
interval $( 0 ,1 )$, because the singularity ${\tilde f}$ is elliptic.
Therefore, using proposition 5.6 we find  that  the  length  of  the
spectrum of the singularity $f$ is less than 1.

\medskip

\noindent
${\rm 6)} \Rightarrow {\rm  2)}$  (cf.  [Tju1])

\smallskip

\noindent
Follows from proposition 5.5.

\medskip

\noindent
{\bf The proof of  corollary
1.}

\smallskip

     It has been proved that the simple singularities are the same  as
the  elliptic  ones.  It  has  been  also  proved  that  an   elliptic
singularity $f$ is stably equivalent to the singularity of a  function
${\tilde f}$ of two variables which has a spectrum of    length
less than one (by proposition 5.6, the length of the spectrum  is  the
same for all stably equivalent singularities).  Now,  by  propositions
5.7 and 5.8,  the  singularity  ${\tilde  f}$  is  equivalent  to  the
singularity of a quasihomogeneous function and the corollary follows.

\bigskip
\bigskip

\vfil
\eject

\centerline {\bf \S 6. The Proofs of Proposition 2.1,  Theorem  2
and Corollary 2.}

\bigskip

\noindent
{\bf The proof of theorem 2.} Theorem 2 can be obtained
immediately from proposition 2.1 and from the following fact.

\noindent
{\bf Proposition 6.1} (see [AC]). {\it For an  isolated
singularity of a  function of $n$ variables the trace  of  the  monodromy
operator is equal to $(-1)^n$. In particular,  if  $n\equiv  3\  ({\rm
mod} \ 4)$ then the trace is equal to $-1$.}

\medskip

\noindent
{\bf The proof of corollary  2.}  In  the  case  of  an
elliptic singularity the action of the braid group $Br(\mu)$ given  by
formula (2) in \S 1 is the action on the set of the tuples of roots in
a root system and it gives the action of $Br(\mu)$ on the tuples  of
reflections in the Weyl group given by formula (3)  in  \S  1  if  one
considers the reflections in the hyperplanes orthogonal to the roots.

     By virtue of these two actions of the braid group $Br(\mu )$, the
corollary follows immediately from the following algebraic  result  of
P.Deligne.

\noindent
{\bf Proposition 6.2} ([De]). {\it Let $W$  be  a  Weyl
group (of one of the types $A_k, D_k$, \allowbreak $E_6 ,  E_7,  E_8$)
of rank $\mu$. Let the tuples $(s_1,\ldots , s_\mu )$ and $(s_1^\prime
,\ldots , s_\mu ^\prime )$ be any elements of $S$ (for the  definition
of $S$ see the end of \S 1) such that $s_1 \cdot \ldots \cdot s_\mu  =
s_1^\prime \cdot \ldots \cdot s_\mu ^\prime =  c$,  where  $c$  is  a
Coxeter element of $W$. Then these two tuples lie in the same orbit of
the action of $Br(\mu )$.}

\medskip

\noindent
{\bf Remark.} This result of P.Deligne has been generalized for
quasi-Coxeter elements (i.e. the ones like $s$ in the
statement of proposition 2.1) by E.Voigt (see [Voi]).

\medskip

\noindent
{\bf  The  proof  of  proposition  2.1.}  The
proposition can be obtained as a consequence of the algebraic theory
of the structure of conjugation classes in a Weyl group (see [Ca]). We
shall briefly outline the proof.

\noindent
{\bf Definition} (see  [Ca],  \S  2).  Let  $w$  be  an
element of $W$. Then the {\it length} $l(w)$ of $w$ is  by  definition
the smallest number $k$ such that $w = w_1  \cdot  \ldots  \cdot  w_k$
where $w_i$, $i= 1,\ldots ,k$, are reflections in $W$.

\noindent
{\bf Proposition 6.3} (see [Ca],  lemma  2).  {\it  The
length $l(w)$ is the number of eigenvalues of $w$ which are not equal to
1.}

\noindent
{\bf Proposition 6.4} (see [Ca], lemma 3).  {\it  If  a
tuple $(s_1,\ldots , s_\mu )$ lies in $S$,  then  the  length  of  the
element $s_1 \cdot \ldots \cdot s_\mu $ is equal to $\mu$.}

\noindent
So the length of the element $s$ from the hypothesis of proposition 2.1
is equal to $\mu$.

\smallskip

By the results  in [Ca] (see [Ca],  \S  3  and
the corollary after proposition 38)  each  element  $w\in  W$  can  be
represented as $w=w_1  \cdot  w_2$,  where  $w_1$  and  $w_2$  can  be
expressed as  pro\-ducts  of  reflections  corresponding  to  mutually
orthogonal roots. The construction described in [Ca],\S  3,  associates
a graph (a so-called {\it admissible  diagram})  $\Upsilon$  to  any
such a representation of $w$. The number of  nodes  in  $\Upsilon$  is
equal to $l(w)$. The admissible diagram for a Coxeter element can be
chosen as the canonical Dynkin diagram of $W$. Conversely, one easily
obtains

\noindent
{\bf Proposition 6.5.} {\it If an admissible diagram of
$w$ is the canonical Dynkin diagram of  $W$  then  $w$  is  a  Coxeter
element of $W$.}

\noindent
Moreover, the following fact turns out to be true.

\noindent
{\bf Proposition 6.6} (see  [Ca],  lemma  8).  {\it  If
$l(w)$ is equal to the rank  of  $W$  and  if  an  admissible  diagram
$\Upsilon$ of $w$ is a tree then $\Upsilon$ is  the  canonical  Dynkin
diagram of $W$.}

\smallskip

\noindent
We also need the following assertion.

\noindent
{\bf Proposition 6.7} (see [Ca],  proposition  22,  and
also proposition 6.3 above). {\it Let $w$ be an element  of  $W$  with
length $l(w)=\mu$ and let $\Upsilon$ be an admissible  diagram  of
$w$. Then the trace of $w$ is given by:  $$\hbox{\rm  tr}  \hbox{w}  =
{\rm  number\ of \ bonds \ in} \ \Upsilon - {\rm  number\ of \
nodes \ in} \ \Upsilon .$$}

\noindent
Now we notice that if $l(w) = \mu $ then an  admissible
diagram corresponding to $w$ is  always  connected  --  otherwise  the
group $W$ would not be irreducible. It is well known that  the  number
of nodes in a connected graph is one greater than the number of  bonds
if and only if the graph is a tree. Hence, by virtue  of  propositions
6.5, 6.6, 6.7,  proposition 2.1   follows.

\bigskip
\bigskip

\centerline {\bf \S 7. The Proof of Theorem II.}

\bigskip

We shall give two proofs of theorem II. Firstly we recall some
notions and results we need for both proofs.

By virtue of corollary  1,  we  can  assume  that  a  simple  (an
elliptic) singularity is quasihomogeneous.

In what follows we shall always assume that \fn is the singularity  of  a
quasihomogeneous function  of  degree  1  with    weights   ${\nu}_1,\ldots
,{\nu}_n$.

\smallskip

\noindent
{\bf  Definition  }  (see   [AGV2]).   Let
$F(z,\lambda )$  ,  $\lambda  \in  \Lambda  =  {\bf  C}^\mu  $,  be  a
miniversal deformation of $f$. Then  the  set  ${\Sigma}_f  =  \lbrace
\lambda \in \Lambda \mid {\rm zero \ is\  a\  critical  \  value\  of}
\break F(\cdot ,\lambda ): ({\bf C}^n , 0) \rightarrow  ({\bf  C},  0)
\rbrace $, $\Sigma \subset \Lambda $, is called
the  {\it  bifurcation
diagram } of $f$.

\smallskip

   Now let $W$ be a Weyl group (of one of the types \ADE )  of  rank
$\mu$. There is a natural action  of  $W$  on  the  complexification
${\bf C}^\mu$ of the Euclidean vector space ${\bf R}^\mu $  on which  $W$
originally acted by  reflections.  The  functions  in  ${\bf  C}[x]$
invariant under the action of $W$ form an algebra which is free  with
$\mu$ generators; these generators (which are also called
{\it the basis invariants } of $W$) can be chosen  as   homogeneous
polynomials of degrees
\break
$m_i +1$, $i=1,\ldots ,\mu$, where  $m_1,\ldots
,m_\mu$ are the exponents of the group
$W$ (see [Ch],[Bou]). Therefore the  space  of  orbits  $B=  {\bf
C}^\mu /W$ is isomorphic to ${\bf C}^\mu$. The set  of  all  singular
orbits (i.e. of all ramification points of the ramified covering  $\pi
: {\bf C}^\mu \rightarrow {\bf C}^\mu / W$) is a hypersurface $S(W)
\subset B$ called the {\it swallowtail } of the group $W$.

\noindent
{\bf Proposition 7.1}  (see  [Lo2],  theorem
4.3). {\it Let \fn be a singularity with    multiplicity
$\mu$ and monodromy group isomorphic to a Weyl group $W$  (of  one
of  the  types  \ADE  )  of  rank  $\mu$.  Then  there  exists  a
biholomorphic map ping of pairs: $(\Lambda ,  {\Sigma}_f)  \rightarrow
(B, S(W)  )  $.  Such  a  biholomorphism  maps  the  fixed  origin  of
coordinates in $B$ into the fixed origin of coordinates in $\Lambda$.}

\smallskip

\medskip

\noindent
{\bf Remark.}
The normal forms of the simple (elliptic)
singularities   are not actually used in the proof of proposition  7.1
(see [Lo2]). It can be checked that the  proof  involves
only the following information on the singularity $f$ (in addition to
the information that the monodromy group of the singularity is an
irreducible Weyl group):

\smallskip

\noindent
a)  $f$  is  (stably)  equivalent  to  the
singularity of a quasihomogeneous function (of three variables);

\smallskip

\noindent
b)  the  quasi\-de\-grees  of  the  ba\-sic
mo\-no\-mi\-als of the lo\-cal al\-ge\-bra $Q_f$ co\-in\-cide with the
num\-bers $\displaystyle{{m_i \over \vert h \vert  }-{\tilde{\nu}}  +1
}$, $i=1,\ldots ,\mu$, where $m_1,\ldots ,m_\mu$ and $\vert  h  \vert$
are  the  exponents  and  the  Coxeter  number  of  the  group   $W$
respectively, \break ${\tilde{\nu}} = {\nu}_1 +\ldots +{\nu}_n$.

\smallskip

\noindent
We have proved  both  facts  not  using  the
normal forms: the  first  one  follows  from
theorem  I  and
corollary 1 and the second one follows from theorem 2 and  proposition
5.10.

\smallskip
\smallskip

So simple singularities with isomorphic monodromy groups have isomorphic
bifurcation diagrams. Then there are two ways to complete the proof of
theorem II. The first one is to use directly the result of
K.Wirthm{\"u}ller ([Wi]) about singularities determined by their
discriminants. The second one is to find a way
to reconstruct the local algebra
$Q_f$ from  the space of orbits $B$ and then to
use the theorem of A.N.Shoshitaishvili (see proposition 7.3)
which says that the singularity
of a quasihomogeneous function is uniquely determined by its local algebra.

\bigskip

\noindent
{\bf The first proof.}

We shall recall some notions we are going to use.
Given a singularity
\break
\gn let $X_g$ denote the germ
of the hypersurface $g^{-1}(0)$ at zero. One can consider
a {\it miniversal deformation of the hypersurface germ} $X_g$ (see
[Tju2], [KS]). It is given by the projection ${\bf C}^n \times {\bf C}^l
\mapsto {\bf C}^l: (z,\lambda ) \mapsto \lambda$, restricted on
the hypersurface $G(z, \lambda ) = 0$, $G(z, \lambda)= g(z) +
{\lambda}_1 e_1(z) + \ldots + {\lambda}_l e_l(z)$,
where $e_1(z),\ldots
e_l(z)$ determine a \hbox{{\bf C}{\hbox{\rm -basis}}} of the vector space
$O_n /{\langle g, {\partial g}/{\partial z_1} ,\ldots ,{\partial g}/
{\partial z_n} \rangle } $, $z\in {\bf C}^n$, $\lambda = ({\lambda}_1,\ldots
,{\lambda}_l )\in {\bf C}^l $. The number $l$ is called the {\it Tjurina
number} of $g$.
The set of critical values of the projection
(i.e. the set of $\lambda \in {\bf C}^l$ such that the variety $G(\cdot
,\lambda) =0$ is singular) forms a hypersurface in ${\bf C}^l$ (more
precisely, one should consider a germ of the hypersurface at $0\in {\bf C}^l$)
called {\it the discriminant of the deformation}.

\smallskip
\smallskip

\noindent
{\bf Proposition 7.2} (the part of the result in [Wi]). {\it
Let $X_1$ and $X_2$ be analytic germs of hypersurfaces with isolated
singularities at $0\in {\bf C}^n$. If the discriminants of some
miniversal deformations of $X_1$ and $X_2$ are isomorphic then
$X_1$ and $X_2$ are isomorphic.}

\smallskip
\smallskip

\noindent
Now we go back to the proof of theorem II. Let \fn and \gn be simple
quasihomogeneous singularities with isomorphic monodromy groups.
Then, by proposition 7.1, $f$ and $g$ have isomorphic bifurcation
diagrams. One can easily see that since $f$ and $g$ are quasihomogeneous
their Tjurina numbers coincide with their multiplicities respectively.
So miniversal
deformations of $f$ and $g$ provide
miniversal deformations of $X_f$
and $X_g$  and the bifurcation diagrams of $f$ and $g$
are isomorphic to their discriminants respectively (since $f$ and $g$
are quasihomogeneous one can actually consider the global discriminant
hypersurface instead of its germ at zero).
Hence, by virtue of proposition 7.2, $X_f$ is isomorphic to $X_g$.
Therefore (see e.g. [Lo3]) there exists a germ of a biholomorphism
of ${\bf C}^n$ mapping $X_f$ onto $X_g$.
It means that $f$ is equivalent to $g h$, where $h$ is an analytic
function invertible at zero. The singularity $g h$ is
semi-quasihomogeneous and since singularity $g$ (which is,
up to a non-zero constant, the principal
part of singularity $g h$) is simple one has,
by the theorem of V.I.Arnold (see [Ar3]), that $g h$ is equivalent
to  $c g$, where $c$ is a non-zero constant. Now again since
$g$ is quasihomogeneous one has that $g$ is equivalent to $c g$.
Thus $f$ is equivalent to $g$ and the theorem is proved.

\bigskip

\noindent
{\bf The second proof.}

Theorem  II  is  the
immediate consequence of the following facts.

\smallskip

\noindent
{\bf Proposition 7.3} (see [Sh]  or  --  for
more general results -- [Ma-Y],[Be]). {\it Let
\break
\fn ,\gn be  singularities of
quasihomogeneous functions.  Then  the  singularities
$f$ and $g$ are equivalent if and only if their local  algebras  $Q_f$
and $Q_g$ are isomorphic (as algebras).}

\smallskip

\noindent
{\bf Lemma.} {\it For a Weyl group  $W$  (of
one of the types \ADE ) such an algebra $A(W)$ can be constructed that
for any simple (elliptic) singularity \fn  with  the  monodromy  group
isomorphic to $W$, the local algebra $Q_f$  is  isomorphic  (as  an
algebra) to $A(W)$.}

\medskip

The proof of the lemma is based on proposition 7.1 and on the V.I.Arnold -
A.B.Givental results on the convolution of invariants of finite groups
generated by reflections. We shall recapitulate  these  results  (for
more details see  [Ar5], [Gi]).

Let $Q_f^\ast$ denote the dual  space  to the  algebra  $Q_f$.
The product of  elements  $p,q$  in  $Q_f$  will be denoted
$p \cdot q$.

\smallskip

\noindent
{\bf Definition} (see [Ar5], [Gi]). A linear
functional $\alpha \in Q_f^\ast$ is said to be {\it ad\-mis\-sible} if
$\alpha$ is not identically equal to zero on the  annihilator  of  the
maximal ideal of $Q_f$ (i.e. on the 1-dimensional ideal generated by
the class of the Hessian of $f$ at zero -- see [AGV1], p.5.11).

\smallskip

A general element  of  $Q_f^\ast$  is  admissible.  If
$\alpha \in Q_f^\ast$ is admissible, then the bilinear  form
$(p,q) \rightarrow \alpha (p\cdot q)$ on $Q_f$  is  nondegenerate.
Let $N_{\alpha } : Q_f \rightarrow Q_f^\ast$ denote  the  operator  of
this form. Let also $\displaystyle{ D={\nu}_1  z_1  {{\partial  }\over
{\partial z_1}}+\ldots + {\nu}_n  z_n  {{\partial  }\over  {\partial
z_n}} }$ be the Euler derivation in the graded local algebra $Q_f$ and
$R=E-D$ be the \hbox{${\bf C}{\hbox{\rm -linear}}$} operator on $Q_f$,
where $E$ is the identity operator.

     By the upper star at an operator  we  shall  denote  the
adjoint ope\-ra\-tor.

     We define a bilinear  operation  $P_{\alpha}  :  Q_f^\ast  \times
Q_f^\ast  \rightarrow  Q_f^\ast$  by   the   formula   \break   $(a,b)
\rightarrow R^\ast N_{\alpha} (N_{\alpha}^{-1} a \cdot N_{\alpha}^{-1}
b)$.

\smallskip

A  function
on the space of orbits
$B$ is called an {\it invariant} of the group $W$. We shall  denote
the tangent space to  $B$  at  $0\in  {\bf  C}^\mu$  by  $T$  and  the
corresponding dual space -- by $T^ \ast$.

     An inner product on ${\bf C}^\mu$ invariant under the  action  of
$W$ provides the isomorphism $i: T^\ast {\bf C}^\mu \rightarrow T_\ast
{\bf C}^\mu$.

     There is a symmetric bilinear operation $\Phi$ on  the  set  of
invariants of $W$ that associates  to  each  pair  $\phi  ,  \psi$  of
invariants the inner product of their Euclidean gradients: $\Phi (\phi
, \psi ) = {\pi}_\ast (i {\pi}^\ast d\phi , i {\pi}^\ast d\psi )$.

     The  invariant  $\Phi  (\phi  ,  \psi  )$  is  called  the   {\it
convolution} of invariants $\phi$ and  $\psi$.  The  operation  $\Phi$
defines a symmetric bilinear operation  ${\Phi}_0  :  T^\ast  \times
T^\ast \rightarrow T^\ast$ by the formula: ${\Phi}_0 (d\phi ,d\psi)  =
d\Phi (\phi , \psi )$. The operation ${\Phi}_0$  is  called  the  {\it
linearized convolution of invariants}.

\smallskip

\noindent
By virtue of proposition 7.1
the space $T^\ast$ can be identified with
the local algebra $Q_f$.

\smallskip


\noindent
{\bf Proposition 7.4} (see [Gi]).
{\it
 Under
any biholomorphism
$(\Lambda , {\Sigma}_f) \rightarrow (B, S(W) ) $ the operation
${\Phi}_0 : T^\ast  \times  T^\ast  \rightarrow  T^\ast$
goes over into the operation
\break
$P_{\alpha}  :  Q_f^\ast  \times
Q_f^\ast  \rightarrow  Q_f^\ast$  for  some  admissible  $\alpha   \in
Q_f^\ast$. Moreover, for any  admissible  $\alpha  \in  Q_f^\ast$  the
operation $P_\alpha$ can be obtained from the operation ${\Phi}_0$  by
means of such a biholomorphism.}


\smallskip

\noindent
{\bf Remark.}

\noindent
1)The normal forms of the simple singularities are not actually
used  in the proof of proposition 7.4 (see [Gi] and the remark
after proposition 7.1).

\noindent
2) Proposition 7.4 was first proved (by  means
of the normal forms) for the simple singularities of the types $A$ and
$D$ (and also for the boundary singularities of the types $B$ and $C$)
in the paper [Ar5].

\medskip

\noindent
{\bf The proof of the lemma} (see [Ar5], \S 9,
Remark 7).

Firstly we shall construct such an algebra $A(W)$. To  do  it  we
apply  the  constructions  and  assertions  mentioned  above  to   the
quasihomogeneous simple (elliptic)
singularity
\fn  and  its  monodromy
group $W$, which is an irreducible Weyl group  of  one  of  the  types
\ADE.

     Let's define the linear operator $w_\beta : T^\ast \to T^\ast$  ,
$\beta \in T^\ast$, by the formula
${  w_\beta  (\cdot  )  =
{\Phi}_0 (\beta ,\cdot )}$.

  As it was said in the beginning
of the paragraph, the basis invariants  of $W$ can be chosen to be
homogeneous polynomials.
Let's  take ${\beta}_0 \in T^\ast$ , ${\beta}_0 = d{\phi}_2 $,
where ${\pi}^\ast {\phi}_2
(z_1,\ldots , z_n) = {z_1}^2 + \ldots +{z_n}^2 $ is the homogeneous
basis
invariant of $W$ of degree 2.
Then, as one easily checks,
the  operator  $w_{{\beta}_0}$  is  in\-ver\-tible.

Now one can
consider  the  family
$A(W)$ of operators $u_\beta = w_{{\beta}_0}^{-1} w_\beta $,
\break
$u_\beta : T^\ast \to T^\ast$,
where  $\beta$  runs  over  the  entire  space
$T^\ast$.
The correspondence $\beta \mapsto u_\beta $ provides $A(W)$ with
the structure of a   \hbox{{$\mu $}{\hbox{\rm -dimensional}}} vector
space. We shall show that
{\it $A(W)$ is the algebra isomorphic to the  local
algebra $Q_f$}.

Indeed,
let's  fix  an
admissible element  $\alpha  \in  Q_f^\ast$  and  let's  fix  such  an
identification  of
$T^\ast$ and $Q_f^\ast$  (see  propositions  7.1,  7.4)
that the operation ${\Phi}_0$ goes over into the operation  $P_\alpha$
under this identification.
Let's define for any  $q\in  Q_f$  the  linear
operator $V_q : Q_f \to Q_f$ by the formula $V_q = M_q R$, where $M_q$
is the operator of multiplication  by  $q$  in  $Q_f$.
Using only the definitions and proposition 7.4  one  checks  (see
[Ar5], \S 9, propositions 3 and 4) that under such  an  identification
$w_\beta = V_q^\ast = R^\ast M_q^\ast $, where ${\beta = N_\alpha q}$.
Notice that since  $ w_{{\beta}_0}$ and $R$ are invertible
$q_0 = N_{\alpha}^{-1} {\beta}_0$ is an invertible element
of $Q_f$.

Now  ${u_\beta = w_{{\beta}_0}^{-1}
w_\beta =}
(M_{q_0}^\ast)^{-1}  (R^\ast)^{-1}  R^\ast  M_q^\ast  =
(M_{q_0}^\ast)^{-1}  M_q^\ast $.
If  $\beta$  runs  over   the   entire space
$T^\ast$, then ${q=N_\alpha^{-1} \beta }$ runs over  the  entire  space
$Q_f$. The operators $ M_q^\ast$
($q$  runs
over the entire space  $Q_f$)  form (with respect to the operator
multiplication)
an  algebra  isomorphic  to ${\widetilde Q}_f$,
where ${\widetilde Q}_f$ is the algebra obtained by
introducing on the vector space $Q_f$ a new
multiplication
opera\-tion
\nobreak $\star$: $a \star b =
q_0^{-1} \cdot a \cdot b$, $a,b \in Q_f$.
One easily checks
that ${\widetilde Q}_f$ is isomorphic  to
$Q_f$ as an algebra.
Hence  $A(W)$ is also an algebra
and it is isomorphic (as an algebra)  to  $Q_f$.

\bigskip
\bigskip

\centerline {\bf Appendix 1. The Normal Forms.}

\bigskip

     Using the information  on  the  simple  (elliptic)  singularities
obtained so far  one  rather  easily  finds  their  normal  forms  and
determines the types of the corresponding Weyl groups.

\noindent
{\bf Theorem } (cf. [Ar1] and [AGV1], \S\S 11,13 ).

\noindent
{\it i) Any simple  (elliptic)  singularity  is  stably
equivalent to one of the following singularities:

     1) $f(x,y) = x^{k+1} + y^2 $, $k\ge 1$;

     2) $f(x,y) = x^2 y + y^{k-1}$, $k\ge 4$;

     3) $f(x,y) = x^3 + y^4$;

     4) $f(x,y) = x^3 + x y^3$;

     5) $f(x,y) = x^3 + y^5$.

\smallskip

\noindent
ii) The Weyl groups  cor\-res\-pon\-ding  to
the sin\-gu\-la\-ri\-ties 1)-5) are of the types $A_k , D_k, E_6,  E_7
, E_8$ res\-pec\-ti\-ve\-ly. In particular,  all  singularities  1)-5)
are mutually not equivalent.}

\medskip

\noindent {\bf The proof.}

\noindent
i) As we know from corollary 1, any  simple  (elliptic)
singularity  is  stably  equi\-va\-lent  to  the  singularity   of   a
quasihomogeneous function of two variables. By virtue  of  proposition
5.9, the diagonal to which all the exponents of the monomials  contained
in that quasihomogeneous  function  belong  lies  below  the  point
$(2,2)$ (in the plane of exponents). So to obtain the normal forms one
sorts out all such lines and uses the following very  easy  assertions
as well as proposition 5.1 ("The Morse lemma with
parameters").

\noindent
{\bf Lemma 1.} {\it Any singularity of  a  function  of
one variable is equivalent to
\break
$f(x)=x^k$ for some $k\ge 2$.}


\noindent
{\bf Lemma 2} (see [AGV1], \S 11.2).  {\it
A cubic form of two variables can  be  reduced  by  a  {\bf  C}-linear
transformation to one of the forms:
     (1) $ x^2 y + y^3 $, (2) $ x^2 y $ ,
\break
(3) $ x^3 $ , (4)  $0
$ .}

\noindent
{\bf Lemma 3} (see [AGV1], \S 11.2). {\it A  polynomial
$A x^2 y + B x y^{k+1} + C y^{2k+1}$, $A\not= 0$, can be reduced by  a
linear transformation to the same form with $B=0$.}

\medskip

\noindent
(ii) We can easily compute  the  spectra  of
the singularities 1)-5) and hence, by means of the  proposition  5.10,
the order and the eigenvalues  of  the  mo\-no\-dro\-my  operator.  By
theorem 2, this enables us to compute  the  exponents  and  the  Coxeter
numbers of
the corresponding Weyl groups. Then we compare these numbers with
ones in the tables (see e.g [Bou]). A Weyl group can be  unambiguously
determined by such data and the theorem follows.

\bigskip
\bigskip

\vfil
\eject

\centerline {\bf Appendix 2. The proof of theorem 1'.}

\bigskip

For the terminology used in this proof see \S 1 and also
[Kou] or [AGV1], [AGV2].

     Let $\Delta$ denote the Newton diagram of $f$. It is  known  (see
[Ar3] or [AGV1]) that one can always  choose  such  a  monomial  basis
$e_1, \ldots , e_\mu$ of $Q_f$  over  ${\bf  C}$  (a  so  called  {\it
regular} basis -- see [Ar3]),  that  for  any  number  $D$  the  basic
monomials of Newton degree $D$ are linearly independent modulo the sum
of  the  gradient  ideal  $I_f$  and  the  space  of  functions  (more
precisely, elements of $O_n$) of  Newton order  greater  than  $D$.
(The Newton degree and the Newton order are defined by  means  of  the
diagram $\Delta$).

     Let's take a miniversal deformation $F$ of  the  singularity  $f$
defined  by  such  a   choice   of   $e_1,   \ldots   ,   e_\mu$,   so
$F(z,\lambda)=f(z)+{{\lambda}_1}e_1+\ldots     +{{\lambda}_\mu}e_\mu$,
\break $\lambda =({\lambda}_1,\ldots ,{\lambda}_\mu) \in {\bf C}^\mu$,
$z\in \Cn$. Let $\[f]$ be now the class of $f$ in $Q_f$ considered  as
a linear combination of monomials
$e_1, \ldots , e_\mu$. One sees that for small $t\in {\bf C}$ the
Newton diagram of the
function $f_t = f+ t \[f]$ coincides with $\Delta$.

     The set $A$ of the principal parts that  are  nondegenerate  with
respect to $\Delta$ is open  in  the  space  of  all  principal  parts
corresponding to  the  diagram  $\Delta$  (see  [Kou]).  Since  for  all
functions with the Newton diagram $\Delta$ and with  principal  part
belonging to the set $A$ the multiplicity of the
critical point $0\in \Cn$
is the same (see [Kou]), the linear 1-parameter deformation $f_t = f+  t
\[f]$ of the singularity $f$ is a $\mu = {\rm const}$ one  (for  small
$t$). The theorem is proved.


\vfil
\eject

\centerline          {\bf          References}

\noindent
\halign{#\hfil&\hfil\quad\quad\vtop{\parindent=0pt\rightskip=1.9cm
\strut#}\cr

     [AC] & A'Campo, N., Le  nombre  de  Lefschetz  d'une  monodromie,
Indag. Math. {\bf 76}:2 (1973), 113-118.
\cr \cr

     [Ar1] & Arnold, V.I., Normal forms for functions near  degenerate
critical points, the Weyl groups $A_k$, $D_k$, $E_k$,  and  Lagrangian
singularities, Funct. Anal. and its Appl. {\bf 6} (1972), 254-272.
\cr \cr

     [Ar2]& Arnold, V.I., Remarks on the stationary phase  method  and
the Coxeter numbers, Russ. Math. Surveys {\bf 28}:5 (1973), 19-48.
\cr \cr

     [Ar3]&  Arnold,  V.I.,  Normal  forms   of   functions   in   the
neighbourhood of degenerate critical points, Russ. Math. Surveys  {\bf
29}:2 (1974), \allowbreak {10-50}.
\cr \cr

     [Ar4]& Arnold, V.I., Critical  points  of  smooth  functions  and
their normal forms, Russ. Math. Surveys {\bf 30}:5 (1975),  1-75.
\cr \cr

     [Ar5]& Arnold, V.I., Indices of singular points of 1-forms  on  a
manifold with  boundary,  the  convolution  of  invariants  of  groups
generated by reflections,  and  the  singular  projections  of  smooth
surfaces, Russ. Math. Surveys {\bf 34}:2 (1979), 1-42.
\cr \cr

     [Ar6]& Arnold, V.I., Catastrophe Theory, Springer, Berlin,  1984.
\cr \cr

     [AGV1]&  Arnold,  V.I.,  Gusein-Zade,  S.M.,   Varchenko,   A.N.,
Sin\-gu\-la\-ri\-ties   of   Diffe\-ren\-tial   Maps,   Vo\-lume    I,
Birkh{\"a}user, Boston-Basel-Stuttgart, 1985.
\cr \cr

     [AGV2]&  Arnold,  V.I.,  Gusein-Zade,  S.M.,   Varchenko,   A.N.,
Sin\-gu\-la\-ri\-ties   of   Diffe\-ren\-tial   Maps,   Vo\-lume   II,
Birkh{\"a}user, Boston-Basel-Stuttgart, 1987.
\cr \cr

     [AGLV]& Arnold, V.I., Goryunov, V.V., Lyashko,  O.V.,  Vassil'ev,
V.A., Singularity Theory I, Encyclopaedia of Math. Sciences,  Vol.  6,
Springer-Verlag, Berlin-New York, 1988.
\cr \cr

     [Be]& Benson, M., Analytic equivalence of  isolated  hypersurface
sin\-gu\-la\-ri\-ties   defined   by   homogeneous   polynomials,   in
"Singularities. Proc. of Symposia in Pure Mathematics. Part 1.", 1983,
111-118.
\cr \cr

     [Bou]& Bourbaki, N., Groupes et alg{\'e}bres de Lie, Ch. IV,V,VI,
Hermann, Paris, 1968.
\cr \cr

    [Br1]& Brieskorn, E., Singular elements of semi-simple algebraic
groups, in "Actes Congr{\'e}s Intern. Math. Nice, 1970", {\bf 2}, Paris,
1971, 279-284.
\cr \cr

    [Br2]&   Brieskorn,   E.,   Die    Monodromie    der    isolierten
Singularit{\"a}ten von Hyper\-fl{\"a}\-chen, Manuscripta Math. {\bf 2}
(1970), 103-161.
\cr \cr

     [Ca]& Carter, R.W., Conjugacy classes in the Weyl group,  Compos.
\allowbreak Math. {\bf 25}:1 (1972), 1-59.
\cr \cr

     [Ch]& Chevalley, C., Invariants of  finite  groups  generated  by
reflections, Amer. J. Math. {\bf 77} (1955), 778-782.
\cr \cr

     [De]& Deligne,  P.,  The  letter  to  E.J.N.Looijenga,  9.3.1974.
Reprinted  in:  
Milnorgitters   einer   einfach    ellipti\-schen    Singularit{\"a}t,
Diplomarbeit, Bonn, 1983, pp. 102-111.
\cr \cr

     [Du]& Durfee, A., Fifteen  chracterizations  of  rational  double
points and simple critical points, Enseignement Math. {\bf 25} (1979),
131-163.
\cr \cr

     [Ga1]&  Gabrielov,  A.M.,  Intersection  matrices   for   certain
singularities, Funct. Anal. and its Appl. {\bf 7} (1973), 182-193.
\cr \cr

     [Ga2]&  Gabrielov,  A.M.,  Bifurcations,  Dynkin   diagrams   and
modality of isolated singularieties, Funct. Anal. and its  Appl.  {\bf
8} (1974), 94-98.
\cr \cr

     [Gi]&  Givental',  A.B.,  Convolution  of  invariants  of  groups
generated by reflections and connected with  simple  singularities  of
functions, Funct. Anal. and its Appl. {\bf 14} (1980), 81-89.
\cr \cr

    [Gr]& Greuel, G.-M., Deformation und Klassifikation von
Singularit{\"a}ten und Moduln, in "Jubil{\"a}um Stagung 100 Jahre DMV
(1990, Bremen)", edited by W.-D.Geyer, Teubner, Stuttgart, 1992, 177-238.
\cr \cr

    [Gu]&   Gusein-Zade,   S.M.,   Monodromy   groups   of   isolated
singularities of hypersurfaces, Russ. Math. Surveys {\bf 32}:2 (1977),
23-69.
\cr \cr

     [Kou]& Kouchnirenko, A.G., Poly{\`e}dres de Newton  et  nombres  de
Milnor, Invent. Math. {\bf 32} (1976), 1-31.
\cr \cr

    [KS]& Kas, A., Schlessinger, M., On the versal deformation of a complex
space with an isolated singularity, Math. Ann. {\bf 196} (1972), 23-29.
\cr \cr

     [Lam]&     Lamotke,     K.,     Die      Homologie      isolierte
Sin\-gu\-la\-ri\-t{\"a}\-ten,  Math.  Zeit.,  {\bf  143}   \allowbreak
{(1975),} 27-44.
\cr \cr

     [Lo1]& Looijenga,  E.J.N.,  The  complement  of  the  bifurcation
variety of a  simple  singularity,  Invent.  Math.  {\bf  23}  (1974),
105-116.
\cr \cr

     [Lo2]&  Looijenga,  E.J.N.,   A   period   mapping   of   certain
semiuniversal deformations, Compos. Math. {\bf 30}:3 (1975),  299-216.
\cr \cr

     [Lo3]&  Looijenga,  E.J.N., Isolated singular points on complete
intersections, London Math. Soc. Lect. Notes Series 77, Cambridge Univ.
Press, Cambridge, 1984.
\cr \cr

     [Ma-Y] & Mather, J.N., Yau, S.S.-T., Classification  of  isolated
hypersurface singularities by their  moduli  algebras,  Invent.  Math.
{\bf 69} (1982), 243-251.
\cr \cr

     [Mi]& Milnor,  J.,  Singular  points  of  complex  hypersurfaces,
Princeton Univ. Press, Providence, 1968.
\cr \cr

     [Sa1]& Saito, K., Quazihomogene isolierte Singularit{\"a}ten  von
Hyper\-fl{\"a}\-chen, Invent. Math. {\bf 14} (1971), 123-142.
\cr \cr

     [Sa2]& Saito, K., A characterization of the intersection form  of
a Milnor fiber for a function with an isolated critical  point,  Proc.
Japan. Acad. Science {\bf 58}:9 (1982), 72-81.
\cr \cr

     [Sa3]&  Saito,  K.,  The  zeroes   of   characteristic   function
${\chi}_f$ for the  exponents  of  a  hypersurface  isolated  singular
point, in "Adv. Studies in Pure Math. I.", 1983, 195-217.
\cr \cr

     [Sh]&  Shoshitaishvili,  A.N.,  Fun\-ctions  with  iso\-mor\-phic
Jaco\-bian ideals, {Funct.} \allowbreak {Anal.} {and} its  Appl.  {\bf
10} (1976), 128-132.
\cr \cr

     [Sl]& Slodowy, P.,  Simple  singularities  and  simple  algebraic
groups, \allowbreak Lecture Notes in Math. {\bf 815} (1980).
\cr \cr

     [St1]&  Steenbrink,  J.H.M.,  Mixed  Hodge  structures   on   the
vanishing cohomology, in "Proceedings of Nordic Summer School on  Real
and  Complex  Singularities",  Sijthoff  and  Noordhoff,  Oslo,  1976,
\allowbreak {525-563}.
\cr \cr

     [St2]&   Steenbrink,   J.H.M.,   The   intersection   form    for
quasihomogeneous  singularities,  Compos.  Math.  {\bf   34}   (1977),
211-223.
\cr \cr

     [Tju1]& Tjurina, G.N., On the topological properties  of  isolated
singularities of complex  spaces,  Math.  USSR  -  Izvestiya  {\bf  2}
(1968), 557-571.
\cr \cr

     [Tju2]& Tjurina, G.N., Locally semiuniversal flat deformations
of isolated singularities of complex spaces, Math. USSR - Izvestija {\bf 3}:5
(1970), 967-999.
\cr \cr

     [V1]& Varchenko,  A.N.,  The  asymptotics  of  holomorphic  forms
determine a mixed Hodge  structure,  Soviet  Math.  Dokl.  {\bf  22}:3
(1980), \break {772-775}.
\cr \cr

     [V2]& Varchenko, A.N., An asymptotic  mixed  Hodge  structure  in
va\-ni\-shing cohomologies, Math. USSR - Izvestija {\bf 18}:3  (1982),
469-512.
\cr \cr

     [V3]& Varchenko, A.N., On the  monodromy  operator  in  vanishing
cohomologies and the operator of multiplication by {$f$} in the  local
ring, Soviet Math. Dokl. {\bf 24}:2 (1981), 248-252.
\cr \cr

     [V4]& Varchenko,  A.N,  The  spectrum  and  decompositions  of  a
critical point of a function, Soviet Math. Dokl.  {\bf  27}:3  (1983),
575-579.
\cr \cr

     [V5]& Varchenko, A.N., A lower bound for the codimension  of  the
stratum \hbox{$\mu {\hbox{\rm =const}}$} in terms of the  mixed  Hodge
structure, Moscow Univ. Math. Bull. {\bf 37}:6 (1982), 30-33.
\cr \cr

     [VC]& Varchenko, A.N., Chmutov, S.V., Tangent cone to the stratum
\hbox{$\mu {\hbox{\rm =const}}$}, Moscow Univ. Math. Bull. {\bf  40}:1
(1985), 7-12.
\cr \cr

     [Voi]& Voigt, E., Ausgezeichnete basen von Milnorgittern einfacher
Singularit{\"a}ten, Abh. Math. Sem. Univ.Hamburg {\bf 55} (1985), 183-190.
\cr \cr

     [Wi]& Wirthm{\"u}ller, K., Singularities determined by their
discriminant, Math. Ann. {\bf 252} (1980), 237-245.
\cr \cr }

\bigskip

DEPARTMENT OF MATHEMATICS, STANFORD UNIVERSITY

STANFORD, CALIFORNIA 94305

\smallskip

{\it E-mail address: } entov@math.stanford.edu

\bye